\documentclass[12pt]{article}
\usepackage{amsmath,amsthm,amsfonts}
\usepackage{graphicx}
\usepackage{fullpage}
\usepackage{url}

\newcommand{\binary}{\lbrace 0,1 \rbrace}

\newtheorem{theorem}{Theorem}

\newtheorem{lemma}[theorem]{Lemma}

\title{Infinite words containing squares at every position}
\author{James Currie\thanks{The author is
supported by an NSERC Discovery Grant.} and
Narad Rampersad\thanks{The author is supported by an NSERC Postdoctoral
Fellowship.} \\
Department of Mathematics and Statistics \\
University of Winnipeg \\
515 Portage Avenue \\
Winnipeg, Manitoba R3B 2E9 (Canada) \\
\url{j.currie@uwinnipeg.ca} \\
\url{n.rampersad@uwinnipeg.ca}}

\begin{document}
\date{\today}
\maketitle

\begin{abstract}
Richomme asked the following question:
what is the infimum of the real numbers $\alpha > 2$ such that
there exists an infinite word that avoids $\alpha$-powers but
contains arbitrarily large squares beginning at every position?
We resolve this question in the case of a binary alphabet by showing
that the answer is $\alpha = 7/3$.
\end{abstract}

\section{Introduction}

We consider the following question of Richomme \cite{Ric05}:
what is the infimum of the real numbers $\alpha > 2$ such that
there exists an infinite word that avoids $\alpha$-powers but
contains arbitrarily large squares beginning at every position?
As we shall see, over the binary alphabet, the answer to Richomme's
question is $\alpha = 7/3$.

First we recall some basic definitions.
If $\alpha$ is a rational number, a word $w$ is an $\alpha$-\emph{power}
if there exist words $x$ and $x'$, with $x'$ a prefix of $x$, such that
$w = x^nx'$ and $\alpha = n + |x'|/|x|$.  We refer to $|x|$ as a
\emph{period} of $w$.  An $\alpha^+$-\emph{power}
is a word that is a $\beta$-power for some $\beta > \alpha$.
A word is $\alpha$-\emph{power-free} (resp. $\alpha^+$-\emph{power-free})
if none of its subwords is an $\alpha$-power (resp. $\alpha^+$-power).
A $2$-power is called a \emph{square}; a $2^+$-power is called an
\emph{overlap}.

The motivation for Richomme's question comes from the observation
that there exist aperiodic infinite binary words that contain
arbitrarily large squares starting at every position.
For instance, all Sturmian words have this property
\cite[Proposition~2]{ADQZ01}.

Certain Sturmian words are not only aperiodic, but also avoid $\alpha$-powers
for some real number $\alpha$.  For instance, it is well-known that the
Fibonacci word
\[
{\bf f} = 010010100100101001010 \cdots
\]
contains no $(2+\varphi)$-powers \cite{MP92} (see also \cite{Kri07}),
where $\varphi = (1+\sqrt{5})/2$ is the golden ratio.  By contrast,
the Thue--Morse word
\[
{\bf t} = 011010011001011010010110 \cdots
\]
is overlap-free \cite{Thu12} but does not contain squares beginning
at every position (It is an easy exercise to show that ${\bf t}$ does not
begin with a square). The squares occurring in the Thue--Morse word have been
characterized by Pansiot \cite{Pan81} and Brlek \cite{Brl89},
and the positions at which they occur were studied by
Brown et al.\ \cite{BRSV06}.

Saari \cite{Saa07} also studied infinite words containing squares
(not necessarily arbitrarily large) beginning at every position.
He calls such words \emph{squareful} (though he imposes the additional
condition that the word contain only finitely many distinct minimal squares).

\section{Overlap-free and $7/3$-power-free squares}

We begin by reviewing what is known concerning the overlap-free binary
squares.  Subsequently, we shall generalize this characterization to the
$7/3$-power-free binary squares.

Let $\mu$ denote the Thue--Morse morphism:
i.e., the morphism that maps $0 \to 01$ and $1 \to 10$.
Define sets \[A = \{00,11,010010,101101\}\] and
\[\mathcal{A} = \bigcup_{k \geq 0}\mu^k(A).\]
The set $\mathcal{A}$ is the set of squares appearing in the Thue--Morse word.
Shelton and Soni \cite{SS85} characterized the overlap-free squares
(the result is also attributed to Thue by Berstel \cite{Ber92}),
as being the conjugates of the words in $\mathcal{A}$ (A conjugate of $x$
is a word $y$ such that $x = uv$ and $y = vu$ for some $u,v$).
An immediate application is the following theorem:

\begin{theorem}
\label{ofree}
If ${\bf w}$ is an infinite overlap-free binary word, then there is
a position $i$ such that ${\bf w}$ does not contain a square beginning
at position $i$.
\end{theorem}

\begin{proof}
An easy computer search suffices to verify that any overlap-free word over
$\binary$ of length greater than $36$ must contain the subword $010011$.
Let $i$ denote any position at which $010011$ occurs in ${\bf w}$.
We claim that no square begins at position $i$.
Suppose to the contrary that $xx$ is such a square.
By Shelton and Soni's result, except for the squares that are conjugates of
$A$, every overlap-free square $xx$ has $|x|$ even.  It follows that $xx$
is of the form $0y10y1$, where $y \in \{01,10\}^\ell$ for some $\ell$.
However, this forces $xx$ to be followed by $0$ in ${\bf w}$, so that
we have the overlap $xx0$ as a subword of ${\bf w}$, a contradiction.
\end{proof}

The next theorem generalizes the characterization of Shelton and Soni.

\begin{theorem}
\label{squares73}
The $7/3$-power-free binary squares are the conjugates of the
words in $\mathcal{A}$.
\end{theorem}

We defer the proof of Theorem~\ref{squares73} to Section~\ref{squares_proof}.  However,
we end this section by proving an analogue, for $7/3$-power-free words, of a well-known
``progression lemma'' for overlap-free words \cite{Ber94,Fif80,Kfo88,Kob88,RS84,SS85}
that we shall need later.

The words $\mu^n(0)$ and $\mu^n(1)$, $n \geq 0$, are known as \emph{Morse blocks}.
Note that the reverse of a Morse block is a Morse block.

\begin{lemma}
\label{prog}
Let $w = uvxy$ be a binary $7/3$-power-free word with
$|u| = |v| = |x| = |y| = 2^n$.  If $u$ and $v$ are Morse blocks,
then $x$ is a Morse block.
\end{lemma}

\begin{proof}
The proof is by induction on $n$.  Clearly, the result holds for
$n = 0$.  We have either (1) $w = u'u''u'u''pqrs$ or (2)
$w = u'u''u''u'pqrs$, where $u'$ and $u''$ are distinct Morse blocks of
length $2^{n-1}$ and $|p|=|q|=|r|=|s|=2^{n-1}$.  By induction, $p$, $q$,
and $r$ are also Morse blocks.  We must show that $p \neq q$.  In case~(1),
$pq = u'u'$ creates the $5/2$-power $u'u''u'u''u'$, and
$pq = u''u''$ creates the cube $u''u''u''$.  In case~(2),
$pq = u'u'$ creates the cube $u'u'u'$; and if $pq = u''u''$, then
$r = u'$ creates the $7/3$-power $u'u''u''u'u''u''u'$, and
$r = u''$ creates the cube $u''u''u''$.  Thus, $p \neq q$, as
required.
\end{proof}

\section{Infinite words containing squares at every position}

We begin by showing that the answer to Richomme's question
is at most $7/3$.

\begin{theorem}
There exists an infinite $(7/3)^+$-power-free binary word that contains
arbitrarily large squares beginning at every position.
\end{theorem}

\begin{proof}
We rely on the existence of a $(7/3)^+$-power-free morphism; that is,
a morphism $f$ such that $f(w)$ is $(7/3)^+$-power-free whenever
$w$ is $(7/3)^+$-power-free.  Kolpakov, Kucherov, and
Tarannikov \cite{KKT99} have given an example of such a morphism
$f$:
\begin{eqnarray*}
0 & \to & 011010011001001101001 \\
1 & \to & 100101100100110010110.
\end{eqnarray*}

For convenience we work instead with the following morphism $g$:
\begin{eqnarray*}
0 & \to & 011010011011001101001 \\
1 & \to & 100101100110110010110,
\end{eqnarray*}
which we obtain by setting $g(0) = \overline{f(1)}$ and
$g(1) = \overline{f(0)}$, where the overline denotes
binary complementation.

We now define a sequence of words $(A_n)_{n \geq 0}$ as follows.
Let $A_0 = 0110110$.  For $n \geq 0$ define $A_{n+1} = (011010)^{-1}g(A_n)$,
where the notation $(011010)^{-1}x$ denotes the word obtained by removing
the prefix $011010$ from the word $x$.  Since $g(A_n)$ always begins
with $011010$, this operation is well-defined.  The sequence $(A_n)_{n \geq 0}$
thus begins
\begin{eqnarray*}
A_0 & = & 0110110 \\
A_1 & = & 011011001101001100101100110110010110100101100110110010110 \cdots \\
& \vdots &
\end{eqnarray*}

Observe that since $g$ is $(7/3)^+$-power-free, $A_n$ is also $(7/3)^+$-power-free.
We first show that as $n \to \infty$, $A_n$ converges to an infinite limit
word ${\bf w}$; second, we show that ${\bf w}$ contains arbitrarily long squares
beginning at every position.

We show by induction on $n$ that $A_n$ is a prefix of $A_{n+1}$.  Clearly
this is true for $n = 0$.  Recall that $A_{n+1} = (011010)^{-1}g(A_n)$.
Inductively, $A_{n-1}$ is a prefix of $A_n$, so $g(A_{n-1})$ is a prefix
of $g(A_n)$.  Thus, $(011010)^{-1}g(A_{n-1}) = A_n$ is a prefix of
$(011010)^{-1}g(A_n) = A_{n+1}$, as required.  We conclude
that $A_n$ tends, in the limit, to a $(7/3)^+$-power-free word ${\bf w}$.

To see that ${\bf w}$ contains arbitrarily long squares beginning at
every position, first let $u = 011010$ and observe that for every $n \geq 1$
we have
\[
A_n = u^{-1}g(u)^{-1}[g^2(u)]^{-1} \cdots [g^{n-1}(u)]^{-1} g^n(0110110).
\]
Let
\[
v = g^{n-1}(u) g^{n-2}(u) \cdots g(u) u,
\]
so that $A_n = v^{-1} g^n(0110110)$,
and observe that
\[
|v| = 6 \sum_{j = 0}^{n-1} 21^j = 6 \left( \frac{21^n - 1}{21 - 1} \right) <
21^n = |g^n(0)|.
\]
We see then that $v$ is a prefix of $g^n(0)$.  Write $g^n(0) = vx$,
so that $A_n = x g^n(11) v x g^n(11) v x$.  It follows that for
$0 \leq j < |x| = 21^n - (3/10)(21^n - 1)$, $A_n$, and hence ${\bf w}$,
contains a square of length $6\cdot 21^n$ beginning at position $j$.
Since $n$ may be taken to be arbitrarily large, the result follows.
\end{proof}

The result of the previous theorem is optimal, as we now demonstrate.

\begin{theorem}
\label{arb_lrg_sqrs}
If ${\bf w}$ is an infinite $7/3$-power-free binary word, then there is
a position $i$ such that ${\bf w}$ does not contain arbitrarily large
squares beginning at position $i$.
\end{theorem}

\begin{proof}
As in the proof of Theorem~\ref{ofree}, an easy computer search suffices
to verify that any $7/3$-power-free word over $\binary$ of length greater
than $39$ must contain the subword $010011$.  Let $i$ be a position at which
there is an occurrence of $010011$ in ${\bf w}$.  Suppose that
there is a square $xx$ beginning at position $i$.  By Theorem~\ref{squares73},
$xx$ is a conjugate of a word in $\mathcal{A}$.  In particular, $xx \notin \mathcal{A}$;
that is, $xx$ is a \emph{non-identity conjugate} of a word in $\mathcal{A}$.

Case~1: $xx$ is a conjugate of either $\mu^k(00)$ or $\mu^k(11)$ for some $k$.
Then $xx$ is a conjugate of a word of the form $uvuv$, where $u$ and $v$
are Morse blocks of the same length.  Without loss of generality, we write
$xx  = u''vuvu'$, where $u'u'' = u$ and $u' \neq \epsilon \neq u''$.

Suppose that $yy$ is another square beginning at position $i$.
Suppose further that there are arbitrarily large squares beginning
at position $i$, so that we may choose $|y| > |xx|$.  We see then
that there is an occurrence of $u''vuvu'$ at position $i + |y|$.
Considering this later occurrence of $u''vuvu'$, and observing that
Morse blocks of a given length are uniquely identified by their first
letter (as well as by their last letter), we may apply Lemma~\ref{prog}
to this later occurrence of $u''vuvu'$ to conclude that the $vuv$
of this occurrence is both preceded and followed by the Morse block $u$.  Thus,
${\bf w}$ contains the $5/2$-power $uvuvu$, a contradiction.

Case~2: $xx$ is a conjugate of either $\mu^k(010010)$ or $\mu^k(101101)$
for some $k$.  By a similar argument as in Case~1, we may suppose that
$xx$ has one of the forms $u''vuuvuu'$, $u''vvuvvu'$, or $u''uvuuvu'$,
where $u$ and $v$ are Morse blocks of the same length, $u'u'' = u$,
and $u' \neq \epsilon \neq u''$.

As before, we suppose the existence of a square $yy$ beginning at
position $i$, where $|y| > |xx|$.  Then there is a later occurrence
of $xx$ at position $i+|y|$.  Applying Lemma~\ref{prog} to this later
occurrence of $xx$, we deduce the existence of one of the $7/3$-powers
$uvuuvuu$, $uvvuvvu$, or $uuvuuvu$, a contradiction.

All cases yield a contradiction; we conclude that there does not exist
arbitrarily large squares beginning at position $i$, as required.
\end{proof}

It is possible, however, to have an infinite $7/3$-power-free binary word with
squares beginning at every position; we are only prevented from having
arbitrarily large squares beginning at every position.

\begin{theorem}
\label{everypos73}
There exists an infinite $7/3$-power-free binary word that contains
squares beginning at every position.
\end{theorem}

\begin{proof}
We show that a word constructed by Currie, Rampersad, and Shallit
\cite{CRS06} has the desired property.  The construction is as follows.
We define the following sequence of words: $A_0 = 00$ and
$A_{n+1} = 0\mu^2(A_n)$, $n \geq 0$.  The first few terms in this sequence
are
\begin{eqnarray*}
A_0 & = & 00 \\
A_1 & = & 001100110 \\
A_2 & = & 0011001101001100101100110100110010110 \\
& \vdots &
\end{eqnarray*}

Currie et al.\ showed that as $n \rightarrow \infty$,
this sequence converges to an infinite word $\mathbf{a}$,
and further, ${\bf a}$ is $7/3$-power-free.  We show that
${\bf a}$ contains squares beginning at every position.
We claim that for $n \geq 0$, ${\bf a}$ contains a word
of the form $xxx'$ at position $(4^n-1)/3$, where
$|x| = 4^{n+1}$ and $x'$ is a prefix of $x$ of
length $4^n$.  Observe that for $n \geq 0$, by the construction of $A_{n+1}$,
we have
\[
A_{n+1} = 0 \mu^2(0) \mu^4(0) \cdots \mu^{2(n-1)}(0) \mu^{2n}(A_1).
\]
However, $\mu^{2n}(A_1) = xxx'$, where $x = \mu^{2n}(0011)$
has length $4^{n+1}$; $x' = \mu^{2n}(0)$ is a prefix of $x$ of
length $4^n$; and $\mu^{2n}(A_1)$ occurs at position
\[
\sum_{i=0}^{n-1} 4^i = \frac{4^n-1}{3},
\]
as claimed.  It follows that for
$i \in [(4^n-1)/3,(4^{n+1}-1)/3-1]$, ${\bf a}$ contains a square of length
$2\cdot 4^{n+1}$ at position $i$.  This completes the proof.
\end{proof}

Although the word constructed in the proof of Theorem~\ref{everypos73}
contains squares beginning at every position, it is not \emph{squareful}
in the sense of Saari \cite{Saa07}, since it does not contain only finitely
many minimal squares.  For a squareful word, there exists a constant
$C$ such that at every position there is a square of length at most $C$.
To see that this does not hold for the word ${\bf a}$ constructed
above, we note that by the factorization theorem of
Karhum\"aki and Shallit \cite{KS04} (Theorem~\ref{fact} below),
any infinite $7/3$-power-free binary word contains
occurrences of $\mu^n(0)$ for arbitrarily large $n$.  However,
$\mu^n(0)$ is a prefix of the Thue--Morse word, and we have already noted
in the introduction that the Thue--Morse word does not begin with a square.
Thus there cannot exist a constant $C$ bounding the length of a minimal
square in ${\bf a}$, so ${\bf a}$ is not squareful.  In general,
no infinite $7/3$-power-free binary word can be squareful.

The result of Theorem~\ref{everypos73} can be strengthened by applying
a more general construction of Currie et al.\ \cite{CRS06}.

\begin{theorem}
For every real number $\alpha > 2$, there exists an infinite $\alpha$-power-free
binary word that contains squares beginning at every position.
\end{theorem}

\begin{proof}
Since Theorem~\ref{everypos73} establishes the result for
$\alpha \geq 7/3$, we only consider $\alpha < 7/3$.
We recall the following construction of Currie et
al.\ \cite[Theorem~14]{CRS06}.
Let $s \geq 3$ and $t \geq 5$ be integers such that
$2 < 3-t/2^s < \alpha$, and such that the word obtained by removing the
prefix of length $t$ from $\mu^s(0)$ begins with $00$.  Let
$\beta = 3-t/2^s$.

We construct sequences of words $A_n$, $B_n$ and $C_n$. Define
$C_0 = 00$ and let $u$ be the prefix of length $t$ of $\mu^s(0)$.
For each $n\geq 0$:

\begin{enumerate}
\item Let $A_n = 0C_n$.
\item Let $B_n = \mu^{s}(A_n)$.
\item Let $C_{n+1} = u^{-1}B_n$.
\end{enumerate}

Currie et al. showed that the $C_n$'s converge to an infinite word
${\bf w}$ that is $\beta^+$-power-free, and hence, $\alpha$-power-free.
For $n \geq 1$, ${\bf w}$ begins with a prefix of the form
\[
u^{-1} \mu^s(0) [\mu^s(u)]^{-1} \mu^{2s}(0) \cdots
[\mu^{ns}(u)]^{-1} \mu^{(n+1)s}(0) \mu^{(n+1)s}(00).
\]
Let $u' = u^{-1}\mu^s(0)$.  Thus, for $n \geq 1$,
${\bf w}$ contains the word $\mu^{ns}(u')\mu^{(n+1)s}(00)$
at position
\[
F_n = 2^s \sum_{i=0}^{n-1} 2^{is} - t \sum_{i=0}^{n-1} 2^{is} =
(2^s - t)\left[\frac{2^{ns} - 1}{2^s - 1}\right],
\]
where $\mu^{ns}(u')$ is a suffix of $\mu^{(n+1)s}(0)$.
Letting $|u'| = t'$ and defining $G_n = |\mu^{ns}(u')| =
t'\cdot 2^{ns}$, we have
\[
F_n = t' \left[\frac{2^{ns} - 1}{2^s - 1}\right] < G_n.
\]
Since $u'$ is a suffix of $\mu^s(0)$, and since ${\bf w}$ begins
with $u'\mu^s(00)$, we see that for $j \in [0,t'-1]$, every subword
of ${\bf w}$ of length $2^{s+1}$ starting at position $j$ is a square.
Similarly, for $n \geq 1$, we see that there is a
square of length $2^{(n+1)s+1}$ starting at position $j$ for
every $j \in [F_n,F_n+G_n-1]$.  Since $F_n < G_n$, there
is thus a square at every position of ${\bf w}$.
\end{proof}

\section{Proof of Theorem~\ref{squares73}}
\label{squares_proof}

 In this section we give the proof of Theorem~\ref{squares73}.
 We begin with some lemmas, but first we recall the
 factorization theorem of Karhum\"aki and Shallit
 \cite{KS04}.

\begin{theorem}[Karhum\"aki and Shallit]
\label{fact}
Let $x\in\binary^*$ be $\alpha$-power-free, $2 < \alpha \leq 7/3$.
Then there exist $u,v\in\{\epsilon,0,1,00,11\}$ and an
$\alpha$-power-free $y\in\binary^*$ such that $x=u\mu(y)v$.
\end{theorem}

\begin{lemma}
\label{oddsq}
Let $xx\in\binary^*$ be $7/3$-power-free.  If $xx = \mu(y)$,
then $|y|$ is even.  Consequently, $y$ is a square.
\end{lemma}

\begin{proof}
Suppose to the contrary that $|y| = |x|$ is odd.  By an exhaustive enumeration
one verifies that $|x| \geq 5$.  But then $xx$ contains two
occurrences of $00$ (or $11$) in positions of different parities,
which is impossible.
\end{proof}

The next lemma is a version of Theorem~\ref{fact} specifically
applicable to squares.

\begin{lemma}
\label{squarefact}
Let $xx\in\binary^*$ be $7/3$-power-free.  If $|xx| > 8$, then either
\begin{itemize}
\item[(a)] $xx = \mu(y)$, where $y\in\binary^*$; or
\item[(b)] $xx = \overline{a}\mu(y)a$, where $a\in\binary$ and
$y\in\binary^*$.
\end{itemize}
\end{lemma}

\begin{proof}
Applying Theorem~\ref{fact}, we write $xx = u\mu(y)v$.
We first show that $|u| = |v| \leq 1$.  Suppose that
$u = 00$.  Then $xx$ begins with one of the words
$000$, $00100$, or $001010$.  The first
and third words contain a $7/3$-power, a contradiction.  The
second word, $00100$, cannot occur later in $xx$, as that
would also imply the existence of a $7/3$-power.  We conclude
$u \neq 00$, and similarly, $u \neq 11$.  A similar
argument also holds for $v$.  Since $|xx|$ is even, we must
therefore have $|u| = |v|$, as required.

If $|u| = |v| = 0$, then we have established (a).
If $|u| = |v| = 1$, it remains to show that $u = \overline{v}$.
If $u = v$, then we have $xx = u \mu(y) u$.  Since $x$ begins and
ends with $u$, we have $x = u \mu(y') = \mu(y'') u$, where $y=y'y''$.
Let $z' = \mu(y')$ and $z'' = \mu(y'')$.  Then $z''$ begins with $u$
and hence with $u \overline{u}$.  This implies that $z'$ begins with
$\overline{u}$ and hence with $\overline{u} u$.  This in turn
implies that $z''$ begins with $u \overline{u} u \overline{u}$,
and hence that $z'$ begins with $\overline{u} u \overline{u} u$.
Now we see that $x$ begins with the $5/2$-power
$u \overline{u} u \overline{u} u$, which is a contradiction.
We conclude that $u \neq v$, as required.
\end{proof}

We are now ready to prove Theorem~\ref{squares73}.

\begin{proof}[Proof of Theorem~\ref{squares73}]
Let $xx$ be a minimal $7/3$-power-free square that is not
a conjugate of a word in $\mathcal{A}$.  That $|xx| > 8$
is easily verified computationally.  Applying Lemma~\ref{squarefact}
leads to two cases.

Case 1: $xx = \mu(y)$.  By Lemma~\ref{oddsq}, $y$ is a square.
Furthermore, $y$ is not a conjugate of a word in $\mathcal{A}$,
contradicting the minimality of $xx$.

Case 2: $xx = \overline{a}\mu(y)a$.  Then
$a\overline{a}\mu(y) = \mu(ay)$ is also a square $zz$.
We show that $zz$ is $7/3$-power-free, and consequently,
by Lemma~\ref{oddsq}, that $ay$ is a $7/3$-power-free square,
contradicting the minimality of $xx$.

Suppose to the contrary that $zz$ contains a $7/3$-power
$s = rrr'$, where $r'$ is a prefix of $r$ and $|r'|/|r| \geq 1/3$.
The word $s$ must occur at the beginning of $zz$ and
we must have $|s| > |x|$; otherwise, $xx$ would contain an occurrence
of $s$, contradicting the assumption that $xx$ is $7/3$-power-free.
We have four cases, depending on the relative sizes of $|r|$ and $|z|$,
as illustrated in Figures~\ref{fig1}--\ref{fig4}.

\begin{figure}[p!]
\begin{center}
\includegraphics[width=5in]{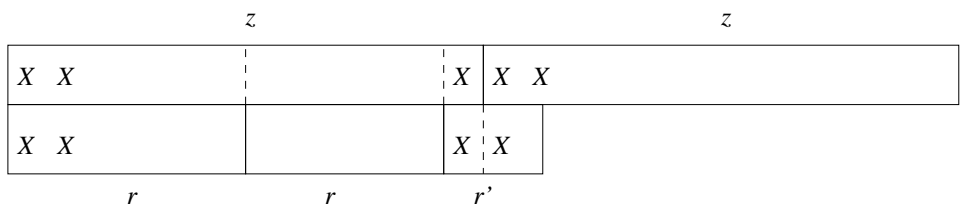}
\end{center}
\caption{The case where $2|r| < |z| \leq 2|r| + |r'|/2$\label{fig1}}
\end{figure}

\begin{figure}[p!]
\begin{center}
\includegraphics[width=5in]{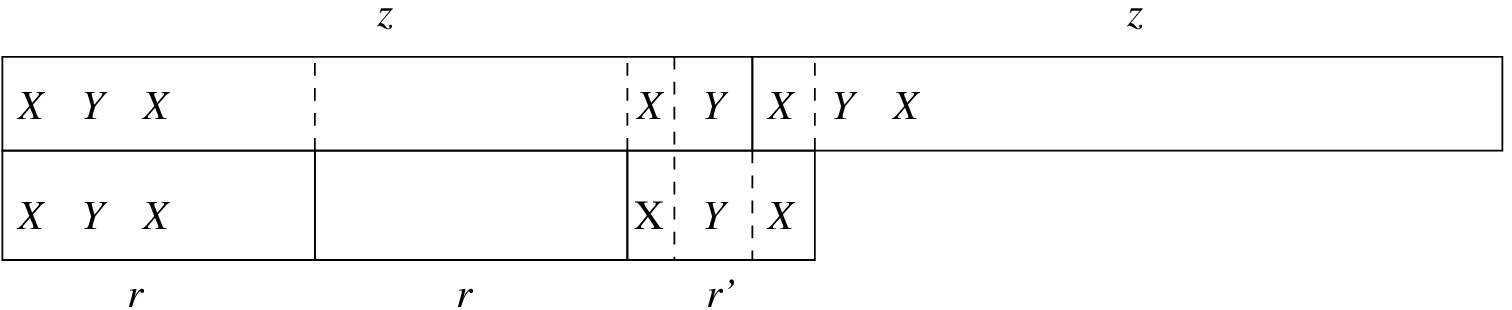}
\end{center}
\caption{The case where $|z| > 2|r| + |r'|/2$\label{fig2}}
\end{figure}

\begin{figure}[p!]
\begin{center}
\includegraphics[width=5in]{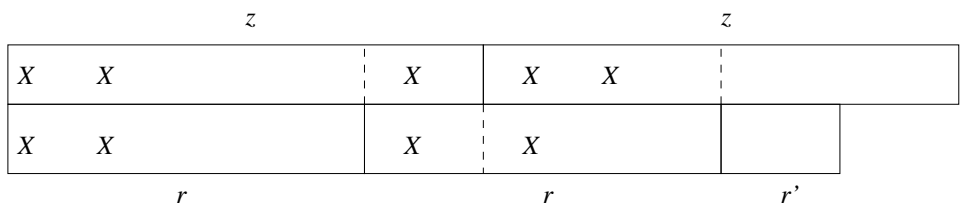}
\end{center}
\caption{The case where $|z| \leq 3/2\cdot |r|$\label{fig3}}
\end{figure}

\begin{figure}[p!]
\begin{center}
\includegraphics[width=5in]{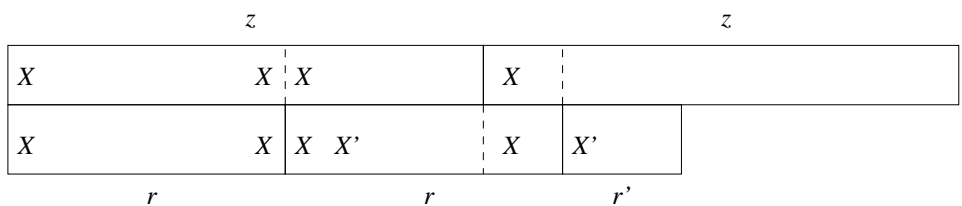}
\end{center}
\caption{The case where $3/2\cdot |r| < |z| < 2|r|$\label{fig4}}
\end{figure}

By analyzing the overlaps between $zz$ and $rrr'$, denoted $X$ in the
figures, we derive a contradiction in each case.

Case~2a: $2|r| < |z| \leq 2|r| + |r'|/2$ (Figure~\ref{fig1}).
In this case, $r'$ has a prefix $X$ that is also a suffix of $z$.  Then $X$ is
also a prefix of $z$ and $XX$ is a prefix of $r'$.  Consequently,
$XX$ is a prefix of $z$ and $xx$ contains the cube $XXX$.  This is
a contradiction.

Case~2b: $|z| > 2|r| + |r'|/2$ (Figure~\ref{fig2}).
In this case, $r'$ is of the form $XYX$, where $XY$ is a suffix of $z$
and $X$ is a prefix of $z$.  Since $r' = XYX$, then $z$ has $XYX$ as a prefix.
In particular, $xx$ contains a $7/3$-power-free square
$XYXY$.  By the assumed minimality of $xx$, $XYXY$ is a conjugate of
a word in $\mathcal{A}$.  If $XYXY = \mu^k(w)$, where $w$ is a conjugate of
a word in $A$ (not $\mathcal{A}$!), then $XYXY$
can be written either as $A_1A_1$ or as $A_1A_2A_3A_1A_2A_3$,
where the $A_i$'s are all Morse blocks of the same length.  Since $XYXY$
is followed by $X$ in $zz$, and since the Morse blocks of a given length
are uniquely identified by their first letter, by Lemma~\ref{prog}
$XYXY$ is followed by the Morse block $A_1$, creating either a cube
or a $7/3$-power in $xx$, contrary to our assumption.

If $XYXY \neq \mu^k(w)$, where $w$ is a conjugate of a word in $A$,
then we may write $XYXY$ as one of $uBABv$, $uBAABAv$, $uBBABBv$, or
$uABAABv$, where $A$ and $B$ are Morse blocks, $B = \overline{A}$,
and $vu = A$.

If $XYXY = uBABv$, then, since $u$ is a non-empty suffix of $A$ and $v$ is
a non-empty prefix of $A$, and since the Morse blocks of a given length
are uniquely identified by their first letter (as well as by their last
letter), we can apply Lemma~\ref{prog} to conclude that $BAB$
is preceded and followed by $A$.  Thus $xx$ contains the
$5/2$-power $ABABA$, contrary to our assumption.
Similarly, if $XYXY = uBAABAv$, then $BAABA$ is preceded and
followed by $A$, creating the $7/3$-power $ABAABAA$.
The other possibilities for $XYXY$ lead to the existence of
a $7/3$-power in $xx$ by a similar argument.

Case~2c: $|z| \leq 3/2\cdot |r|$ (Figure~\ref{fig3}).
In this case, $r$ has a prefix $X$ that is also a suffix of $z$.  Then $X$ is
also a prefix of $z$ and $XX$ is a prefix of $r$.  Consequently,
$XX$ is a prefix of $z$ and $xx$ contains the cube $XXX$.  This is
a contradiction.

Case~2d: $3/2\cdot |r| < |z| < 2|r|$ (Figure~\ref{fig4}).
In this case, $r$ has a suffix $X$ that is also a prefix of $z$.  Then $X$ is
also a prefix of $r$ and $r'$ begins with some prefix $X'$ of $X$.
Consequently, $rr$, and indeed $xx$, contains a subword $XXX'$
that is at least a $7/3$-power.  This is a contradiction.

Since in all cases we have derived a contradiction by showing that
$xx$ contains a $7/3$-power, we conclude that our assumption
that $zz$ contains a $7/3$-power is false.  Recalling that
$zz = \mu(ay)$ and that $ay$ is necessarily a square, we conclude
that $ay$ is a $7/3$-power-free square, contradicting the minimality
of $xx$.  We conclude that there exists no $7/3$-power-free square
$xx$ that is not a conjugate of a word in $\mathcal{A}$.
\end{proof}

We claim that the constant $7/3$ in Theorem~\ref{squares73} is best possible.
To see this, note that the word
\[
01101001101100101100110100110110010110
\]
is a $(7/3)^+$-power-free square, but is not a conjugate of a word
in $\mathcal{A}$.

\section{Conclusion}

We have only considered words over a binary alphabet.  It remains
to consider whether similar results hold over a larger alphabet.
For instance, does there exist an infinite overlap-free ternary
word that contains squares beginning at every position?
Richomme observes that over any alphabet there cannot exist
an infinite overlap-free word containing infinitely many squares
at every position.  He points out that this follows from the following
result of Ilie \cite[Lemma~2]{Ili07}: In any word, if $vv$ and $uu$ are two
squares at position $i$ and $ww$ is a square at position $i+1$, then either
$|w| = |u|$ or $|w| = |v|$ or $|w| \geq 2|v|$.  An easy consequence of
this result is that in any infinite word, if infinitely many distinct
squares begin at position $i$ and $ww$ is a square beginning at position
$i+1$, then $|w| = |u|$ for some square $uu$ occurring at position $i$,
and hence there is an overlap at position $i$.

\section*{Acknowledgments}

The authors wish to thank Gw\'ena\"el Richomme for suggesting
the problem.  We also thank him for reading an earlier draft
of this paper and providing many helpful comments and suggestions.


\begin{thebibliography}{99}
\bibitem{ADQZ01}
J.-P. Allouche, J. L. Davison, M. Queff\'elec, L. Q. Zamboni,
``Transcendence of Sturmian or morphic continued fractions'',
\textit{J. Number Theory} \textbf{91} (2001), 39--66.

\bibitem{Ber92}
J. Berstel, ``Axel Thue's work on repetitions in words''. In
P. Leroux, C. Reutenauer, eds., \emph{S\'eries formelles et combinatoire
alg\'ebrique}, Publications du LaCIM, pp 65--80, UQAM, 1992.

\bibitem{Ber94}
J. Berstel, ``A rewriting of Fife's theorem about overlap-free words''.
In J. Karhum\"aki, H. Maurer, G. Rozenberg, eds., \emph{Results
and Trends in Theoretical Computer Science}, LNCS 812, pp. 19--29
Springer-Verlag, 1994.

\bibitem{Brl89}
S. Brlek, ``Enumeration of factors in the Thue-Morse word'',
\emph{Discrete Appl. Math.} \textbf{24} (1989), 83--96.

\bibitem{BRSV06}
S. Brown, N. Rampersad, J. Shallit, T. Vasiga, ``Squares and overlaps in
the Thue-Morse sequence and some variants'', \textit{Theor. Inform. Appl.}
\textbf{40} (2006), 473-484.

\bibitem{CRS06}
J. Currie, N. Rampersad, J. Shallit,  ``Binary words containing infinitely
many overlaps'', \textit{Electron. J. Combinatorics} \textbf{13} (2006), \#R82.

\bibitem{Fif80}
E. Fife, ``Binary sequences which contain no $BBb$'', \emph{Trans. Amer.
Math. Soc.} \textbf{261} (1980), 115--136.

\bibitem{Ili07}
L. Ilie, ``A note on the number of squares in a word'',
\textit{Theoret. Comput. Sci.} \textbf{380} (2007), 373--376.

\bibitem{KS04}
J. Karhum\"aki, J. Shallit, ``Polynomial versus exponential growth
in repetition-free binary words'', \emph{J. Combin. Theory Ser. A}
\textbf{104} (2004), 335--347.

\bibitem{Kfo88}
R. Kfoury, ``A linear time algorithm to decide whether a binary word
contains an overlap'', \emph{Theoret. Inform. Appl.} \textbf{22} (1988),
135--145.

\bibitem{Kob88}
Y. Kobayashi, ``Enumeration of irreducible binary words'',
\emph{Discrete Appl. Math.} \textbf{20} (1988), 221--232.

\bibitem{KKT99}
R. Kolpakov, G. Kucherov, Y. Tarannikov, ``On repetition-free binary words
of minimal density'', WORDS (Rouen, 1997), \emph{Theoret. Comput. Sci.}
\textbf{218} (1999), 161--175.

\bibitem{Kri07}
D. Krieger, ``On critical exponents in fixed points of non-erasing morphisms'',
\emph{Theoret. Comput. Sci.} \textbf{376} (2007), 70--88.

\bibitem{MP92}
F. Mignosi, G. Pirillo, ``Repetitions in the Fibonacci infinite word'',
\emph{RAIRO Inform. Th\'eor.} \textbf{26} (1992), 199-–204.

\bibitem{Pan81}
J. J. Pansiot, ``The Morse sequence and iterated morphisms'',
\emph{Inform. Process. Lett.} \textbf{12} (1981), 68--70.

\bibitem{RS84}
A. Restivo, S. Salemi, ``Overlap free words on two symbols''. In M. Nivat,
D. Perrin, eds., \emph{Automata on Infinite Words}, LNCS 192,
pp. 198--206, Springer-Verlag, 1984.

\bibitem{Ric05}
G. Richomme. Personal communication, 2005.

\bibitem{Saa07}
K. Saari, ``Everywhere $\alpha$-repetitive sequences and Sturmian words''.
In \textit{Proc. CSR 2007}, LNCS 4649, pp. 362–-372, Springer-Verlag, 2007.

\bibitem{SS85}
R. Shelton, R. Soni, ``Chains and fixing blocks in irreducible
sequences'', \emph{Discrete Math.} \textbf{54} (1985), 93--99.

\bibitem{Thu12}
A. Thue, ``\"{U}ber die gegenseitige Lage gleicher Teile gewisser
Zeichenreihen'', \textit{Kra. Vidensk. Selsk. Skrifter. I. Math. Nat. Kl.}
\textbf{1} (1912), 1--67.
\end{thebibliography}
\end{document}